\documentstyle[12pt]{article}
\newcommand{\Sub}{\mathop{\rm Sub}\limits}

\newcommand{\const}{\mathop{\rm const}\limits}

\newcommand{\argmax}{\mathop{\rm argmax }\limits}

\begin{document}
 \begin{center}

{\bf  EXACT VALUE FOR SUBGAUSSIAN NORM OF \\

\vspace{4mm}

CENTERED INDICATOR RANDOM VARIABLE}\\

\vspace{4mm}

{\sc Eugene Ostrovsky, Leonid Sirota}\\

\vspace{3mm}

 Bar-Ilan University,  59200, Ramat Gan, ISRAEL; \\

 e-mail: eugostrovsky@list.ru \\
 e-mail: sirota3@bezeqint.net \\

 \vspace{4mm}

        {\bf Abstract}

\end{center}

\vspace{3mm}

  We calculate the exact subgaussian norm of a centered (shifted) indicator (Bernoulli's) random variable.\par
 Using this result we derive very simple tail estimates for sums of these variables, not necessary to be
 identical distributed, and give some examples to show the exactness of our estimates. \par

\vspace{3mm}


\vspace{3mm}

 \section { Introduction. Notations. Statement of problem}

 Let $  \{\Omega, B, {\bf P}  \} $ be some non-trivial probability space.
 We say that the {\it centered:} $ {\bf E} \xi = 0 $ numerical random variable (r.v.)
 $ \xi = \xi(\omega), \ \omega \in \Omega $ is subgaussian, or equally, belongs to the space $ \Sub(\Omega), $
 if there exists some non-negative constant $ \tau \ge 0 $ such that

$$
\forall \lambda \in R  \ \Rightarrow
{\bf E} \exp(\lambda \xi) \le \exp[ \lambda^2 \ \tau^2 ]. \eqno(1.1).
$$
 The minimal value $ \tau $ satisfying (1.1) is called a  subgaussian  norm
of the variable $ \xi, $ write

 $$
 ||\xi||\Sub = \inf \{ \tau, \ \tau > 0: \ \forall \lambda \in R \ \Rightarrow {\bf E}\exp(\lambda \xi) \le \exp(\lambda^2 \ \tau^2) \}.
 $$

 Evidently,

$$
||\xi||\Sub = \sup_{\lambda \ne 0} \left[ \sqrt{ \ln {\bf E}  \exp ( \lambda \xi)  }/|\lambda| \right].  \eqno(1.2)
$$

 This important notion was introduced by  J.P.Kahane \cite{Kahane1}; V.V.Buldygin and Yu.V.Kozatchenko \cite{Buldygin1} proved
that the set $ \Sub(\Omega) $  relative the norm $  ||\cdot|| $ is complete Banach space which is isomorphic to subspace
consisting only from the centered variables of Orlicz's space over $ (\Omega, B,P)  $ with $ N \ - $ Orlicz-Young function
 $ N(u) = \exp(u^2) - 1 $  \cite{Kozatchenko1}.  \par

   If $ ||\xi||\Sub = \tau \in (0,\infty),  $ then

 $$
 \max [{\bf P}(\xi > x),  {\bf P}(\xi < -x)  ] \le \exp(- x^2/(4 \tau^2)  ), \ x \ge 0; \eqno(1.3)
 $$
 and  the last inequality is in general case non-improvable.  It is sufficient for this to consider the case when
 the r.v. $  \xi  $ has the centered Gaussian non-degenerate distribution.\par

  Conversely, if  $ {\bf E} \xi = 0 $ and if  for some positive finite constant $  K  $

 $$
 \max [{\bf P}(\xi > x),  {\bf P}(\xi < -x)  ] \le \exp(- x^2/K^2  ), \ x \ge 0,
 $$
 then $ \xi \in \Sub(\Omega) $ and $ ||\xi||\Sub < 4 K. $ \par

 The subgaussian norm in the subspace of the centered r.v. is equivalent to the following Grand Lebesgue Space (GLS)
 norm:

 $$
|||\xi||| := \sup_{s \ge 1} \left[ \frac{|\xi|_s}{\sqrt{s}} \right], \hspace{6mm} |\xi|_s =  \left[ {\bf E} |\xi|^s \right]^{1/s}.
 $$

 For the non - centered r.v. $ \xi $  the subgaussian norm may be defined as follows:

 $$
 ||\xi|| \Sub := \left[  \left\{ ||\xi - {\bf E} \xi||\Sub \right\}^2 + ( {\bf E} \xi)^2  \right]^{1/2}.
  $$

 More detail investigation of these spaces see in the monograph \cite{Ostrovsky1}, chapter 1.  \par

 We denote as usually by $  I(A) = I(A; \omega), \ \omega \in \Omega, \ A \in B  $ the indicator function of event $  A. $
 Further, let $  p  $ be arbitrary number from the set $ [0,1]: \ 0 \le p \le 1 $ and let $ A(p) $  be any event such that
 $ {\bf P}(A(p)) = p. $ Denote also  $ \eta_p = I(A(p)) - p; $  the {\it  centering  } of the r.v. $ I(A(p)); $
 then $ {\bf E} \eta_p = 0  $ and

 $$
 {\bf P}(\eta_p = 1 - p) = p; \hspace{6mm}  {\bf P}(\eta_p = -p) = 1 - p.   \eqno(1.4)
 $$

\vspace{3mm}

{\bf  Our aim in this short report is to compute the exact value of the subgaussian norm for the random variable $ \eta_p. $ }

\vspace{3mm}

 We apply the obtained estimate in the third section to the tail computation for sums of independent indicators.\par

\vspace{3mm}

 Applications of these estimates in the non-parametrical statistics may be found in the articles
\cite{Gaivoronsky1}, \cite{Kiefer1}. Another  application is described in \cite{Chen1}. \par

 \vspace{3mm}

 \section{Main result}

 \vspace{3mm}

 Define the following non-negative continuous on the closed segment $ p \in [0,1]  $ function

$$
Q(p) =  \sqrt{ \frac{1 - 2p}{4 \ln(( 1 - p )/p)} },
\eqno(2.0)
$$
so that $ Q(0+0) = Q(1-0) = 0 $  and $ Q^2(1/2) = 1/8 $ (Hospital's rule).  Note also

$$
p \to 0+ \ \Rightarrow Q(p) \sim \frac{0.5}{\sqrt{|\ln p|}}, \hspace{6mm}  p \to 1 - 0 \Rightarrow  Q(p) \sim \frac{0.5}{\sqrt{|\ln(1 - p)|}}.
\eqno(2.1)
$$
 The last circumstance play a very important role in the non - parametrical statistics,  see \cite{Gaivoronsky1}, \cite{Kiefer1}.\par

\vspace{3mm}

 {\bf Theorem 2.1.}

 $$
 ||\eta_p||\Sub = Q(p). \eqno(2.2)
 $$

\vspace{3mm}

{\bf Proof.}  Note that

$$
{\bf E} e^{\lambda \eta_p} = p e^{\lambda (1 - p) } + (1 - p) e^{ - p \lambda }, \lambda \in (-\infty, \infty).
$$

The inequality

$$
p e^{\lambda (1 - p) } + (1 - p) e^{ - p \lambda } \le e^{ Q^2(p) \ \lambda^2  } \eqno(2.3)
$$
is proved in   \cite{Kearns1}; see also \cite{Berend1}.  Therefore

$$
{\bf E} e^{\lambda \eta_p} \le e^{ Q^2(p) \ \lambda^2  }.
$$
 This imply by direct definition of the subgaussian norm that $ ||\eta||\Sub \le Q(p), \ p \in [0,1]. $\par
Let us prove the inverse inequality.   Suppose $ p \in (0,1); $ the extremal cases $ p = 0, \ p = 1 $ are trivial.\par

We denote following  the authors of articles \cite{Berend1}, \cite{Schlemm1}

$$
\lambda_0 = \lambda_0(p) = 2 \log \left[ \frac{1-p}{p} \right] \eqno(2.4)
$$
and deduce after simple calculations

$$
|| \eta_p || \Sub = \sup_{\lambda \ne 0} \left[ \sqrt{ \ln {\bf E}  \exp ( \lambda \eta_p)  }/|\lambda| \right] =
$$

$$
\sup_{\lambda \ne 0} \left[ \sqrt{ \ln \left\{ p e^{\lambda (1 - p)} + (1 - p) e^{-\lambda p} \right\}  }   /|\lambda| \right] \ge
$$

$$
  \sqrt{ \ln \left\{ p e^{\lambda_0 (1 - p)} + (1 - p) e^{-\lambda_0 p} \right\} } /|\lambda_0|  = Q(p), \eqno(2.5)
$$
Q.E.D.\par

\vspace{3mm}

{\bf Remark 2.1.} Let us explain the choice of the value $  \lambda_0 = \lambda_0(p), \ 0 < p < 1. $ In accordance to the
equality (1.2) the optimal value of the parameter $ \lambda $ is following: $ \lambda = \Lambda(p), $ where

$$
\Lambda(p) = \argmax_{\lambda \ne 0} \left[ \sqrt{ \ln {\bf E}  \exp ( \lambda \eta_p)  }/|\lambda| \right]
$$
or equally

$$
\Lambda(p) = \argmax_{\lambda \ne 0} \left\{ \frac{ \ln ( p e^{ ( 1 - p) \lambda } + (1 - p) e^{-\lambda p}  )}{\lambda^2} \right\}. \eqno(2.6)
$$

 Denote following  the authors of articles \cite{Berend1}, \cite{Kearns1}, \cite{Schlemm1}

$$
g(\lambda) = g_p(\lambda) =  \frac{ \ln ( p e^{ ( 1 - p) \lambda } + (1 - p) e^{-\lambda p}  )}{\lambda^2}.
$$
  It is easy to verify that $ g'_{\lambda}(\lambda_0) = 0. $
 It is also proved in the article  \cite{Schlemm1}  that the function $ \lambda  \to g_p(\lambda) $ is unimodal.
 Therefore, we derive taking into account the behavior of the function $ g_p(\lambda) $  at $ \lambda \to \pm \infty $
that $ \Lambda(p) = \lambda_0. $ \par

 \vspace{3mm}

{\bf  Consequence 2.1. } Let $ \nu:  \Omega \to R $ be a centered stepwise (simple) r.v. (measurable function):

$$
\nu = \sum_{j=1}^{m} c(j) [ I(A(p(j))) - p(j)  ], \ m = \const \le \infty,  \ c(j) = \const, \eqno(2.7)
$$
and $ \{ A(p(j)) \} $ are events not necessary to be disjoint or independent. We conclude  using triangle inequality
for the subgaussian norm and the completeness of the space $  \Sub(\Omega) $  in the case when $  m = \infty: $

$$
||\nu||\Sub \le \sum_{j=1}^m |c(j)| Q(p(j)). \eqno(2.8)
$$

\vspace{3mm}

 \section{Tail estimations for sums of independent indicators}

 \vspace{3mm}

 Let $  p(i), \ i = 1,2, \ldots,n  $ be positive numbers such that $ 0 < p(i) < 1,  $ and let $ A(i) $ be {\it independent} events
for which $  {\bf P}(A(i)) = p(i). $  Introduce a sequence of two - values independent random variables
$ \zeta(i) = I(A(i)) - p(i), $ and define its sum

$$
 S(n) := \sum_{i=1}^n \zeta(i).\eqno(3.1)
$$

 It is known \cite{Ostrovsky1}, chapter 1, section 1.6 that

$$
||\sum_{i=1}^n \zeta(i)||\Sub \le \sqrt{ \sum_{i=1}^n (||\zeta(i)||\Sub)^2  }. \eqno(3.2)
$$

 Therefore

$$
||\sum_{i=1}^n \zeta(i)||\Sub \le W(n), \eqno(3.3)
$$
where
$$
W(n) \stackrel{def}{=} \sqrt{ \sum_{i=1}^n (Q(p(i)))^2  }. \eqno(3.4)
$$
 As a consequence: \\

 \vspace{3mm}

 {\bf Proposition 3.1.}

$$
\max [{\bf P}(S(n) > x),  {\bf P}(S(n) < -x) ] \le \exp \left( -x^2/(4 W^2(n))   \right), \ x \ge 0. \eqno(3.5)
$$

\vspace{3mm}

{\bf Example 3.1.}  Assume in addition that $ p(i) = p = \const \in (0,1); $ then the r.v. $ S(n) + n p $ has a (non-degenerate)
Bernoulli distribution. We deduce using inequalities (3.3) and (3.5):

 $$
\sup_n ||S(n)/\sqrt{n}||\Sub \le Q(p).\eqno(3.6)
 $$
As long as

 $$
\sup_n ||S(n)/\sqrt{n}||\Sub \ge ||S(1)|| \Sub =  Q(p), \eqno(3.6a)
 $$
we get:

 $$
\sup_n ||S(n)/\sqrt{n}||\Sub = Q(p).\eqno(3.7)
 $$
 Thus, the estimate  (3.3) is non-improvable;  cf. \cite{Bentkus1}, \cite{Pinelis1}, \cite{Serov1}, \cite{Zubkov1}. \par

\vspace{3mm}

{\bf Example 3.2.} Suppose in addition to the previous example 3.1 that $  p = 1/2 $ (symmetrical case);
 then it follows from Proposition (3.1) alike the famous Hoeffding's inequality

$$
 {\bf P} ( 2 S(n)/\sqrt{n} > x ) \le  e^{-x^2/2},  \ x > 0,
$$

while

$$
\sup_n {\bf P} ( 2 S(n)/\sqrt{n} > x ) \ge  \lim_{n \to \infty} {\bf P} ( 2 S(n)/\sqrt{n} > x ) \ge C x^{-1} \ e^{-x^2/2}, \
x \ge 1.
$$

\vspace{3mm}

 Another generalizations of the equality (2.2), for example, on the Hoeffding's inequality and on the theory of martingales
see in the article of M.Raginsky and I.Sason \cite{Raginsky1}.\par

\end{document}